# Congestion management via increasing integration of electric and thermal energy infrastructures

*Alvaro* Gonzalez-Castellanos[1*], *Priyanko* Guha Thakurta[2], and *Aldo* Bischi[1]

[1]Center for Energy Science and Technology, Skolkovo Institute of Science and Technology, Moscow, Russia
[2]EirGrid plc, Dublin, Ireland

**Abstract.** Congestion caused in the electrical network due to renewable generation can be effectively managed by integrating electric and thermal infrastructures, the latter being represented by large scale District Heating (DH) networks, often fed by large combined heat and power (CHP) plants. The CHP plants could further improve the profit margin of district heating multi-utilities by selling electricity in the power market by adjusting the ratio between generated heat and power. The latter is possible only for certain CHP plants, which allow decoupling the two commodities generation, namely the ones provided by two independent variables (degrees-of-freedom) or by integrating them with thermal energy storage and Power-to-Heat (P2H) units. CHP units can, therefore, help in the congestion management of the electricity network. A detailed mixed-integer linear programming (MILP) optimization model is introduced for solving the network-constrained unit commitment of integrated electric and thermal infrastructures. The developed model contains a detailed characterization of the useful effects of CHP units, i.e., heat and power, as a function of one and two independent variables. A lossless DC flow approximation models the electricity transmission network. The district heating model includes the use of gas boilers, electric boilers, and thermal energy storage. The conducted studies on IEEE 24 bus system highlight the importance of a comprehensive analysis of multi-energy systems to harness the flexibility derived from the joint operation of electric and heat sectors and managing congestion in the electrical network.

## 1 Introduction

The massive deployment of renewable energy sources (RES) vastly reduces greenhouse emissions and operating costs. However, their integration into power systems is constrained by their geographical availability, not necessarily matching demand location. These geographical constraints can result in transmission congestion if the network capacity is not enough to transfer power from RES sites to load centres. Transmission congestion can lead to curtailment of RES, generation rescheduling, or load shedding, thereby leading to economic losses for the system, loads, and both fuel-based and RES generators.

Congestion management (CM) is conventionally performed via technical and market measures [1]. Technical measures include modifying network parameters, such as the creation or enhancement of transmission corridors, line-switching, flexible AC transmission system (FACTS) devices, and management of substations' topology, amongst others. Market measures are characterized by market products such as spinning reserves and demand-side management, or the use of pricing mechanisms like nodal pricing, uplift costs, and price area CM. Even though these measures have proven effective for CM, they result in significant capital investments or economic losses. The cost of any measure used for CM should compensate for the congestion in the system [2].

A measure that would aid in CM and increase network flexibility is the operational integration of the heating and electricity sectors. It can be done by optimal management of CHP plants, whose flexibility can be boosted via heat storage integration and Power-to-Heat (P2H) solutions. In this way, the systems' short-term (up to daily cycles for heat storage) flexibility is increased balancing of RES fluctuations and reducing operational costs [3].

Combined heat and power (CHP) units produce heat and electricity depending on their independent variables, i.e., degrees-of-freedom. The independent variable of so-called one-degree-of-freedom units is their fuel consumption, e.g., a simple cycle gas turbine or a backpressure steam turbine. These units cannot easily decouple thermal and electric power. They increase or decrease jointly both with fuel consumption.

---

* Corresponding author: alvaro.gonzalez@skolkovotech.ru

Other units can be controlled based on two-degrees-of-freedom, in which a second independent variable allows to decouple the heat and electricity generation. An example of this type of unit is an extraction condensing steam turbine whose independent variables are the fuel and the valve opening controlling the ratio of heat and electric generation depending on its profitability. Such generation flexibility could balance RES fluctuations, either directly if they belong to the generation portfolio of the multi-utility or indirectly reacting to electric prices or required RES-balancing services. CHP units can provide a higher thermal output when RES generation is higher and contribute to the electric energy balance when RES decreases. Such flexibility can be obtained from one-degree-of-freedom units via heat storage and P2H.

The conventional practice of operating large-scale CHPs is a "decoupled" one, i.e., the units are primarily following the thermal load. After the thermal working point has been set, an electric economic dispatch is performed. The committed heat generation thus constrains the generated electricity. This approach does not achieve the minimum global cost for the integrated energy systems operation since the optimization of the thermal and electric systems is performed in a sequential rather than integrated way. When accounting for energy storage or P2H, a decoupled dispatch further diverges from the optimal operation [4].

Several works addressed the integration of electric power flow in integrated heat and power models. Among them are methods based on developing non-convex models for the electric and thermal networks [5] and the integration of gas networks [6]. Additional works have focused on developing convex models based on multi-energy virtual power plants, i.e., distributed energy generation operated as one larger plant; and integration of RES [7]. These models aimed to analyze the energy flow between the systems and economics, neglecting the effect of the network constraints on the flexible operation of the units.

Furthermore, in the above-referenced works, the generation plants are modeled as constant efficiency units, strongly simplifying the nonlinearity of their performance curves. An accurate representation of the generation at partial loads becomes essential when considering short-term scenarios and assessing system flexibility [8]. Constant efficiency modeling can result in over and under generation.

Rong *et al.* developed a model for a CHP dispatch with multiple generation and consumption sites based on dynamic programming [9]. In this study, the electric power flow is based on an energy flow model, surplus and shortage, rather than on the actual path that the power would follow due to the electric characteristics of transmission lines. This can lead to overestimation of the amount of power flowing between lines, making it necessary to re-dispatch some of the scheduled generators and incurring compensation costs. Morvaj *et al.* presented an optimization framework for integrating thermal and electric energy systems while introducing electric network constraints [10]. However, this study was performed for electric distribution networks, which cover smaller areas and are not as constrained in their RES location.

The main contribution of this paper is the formulation of a model that allows the assessment of the flexibility impact derived from the integration of electric and thermal infrastructures via large scale CHP plants, thermal storage, and P2H. For this purpose, a mixed-integer linear programming (MILP) model for unit commitment is presented that characterizes in a detailed manner the non-convex performance curves of electric generators, gas boilers, cogeneration units with one and two degrees-of-freedom, thermal energy storage, and power-to-heat units. The P2H units are assumed to be electric boilers in order to be sure that the more common DH networks with pressurized water temperatures above 100ºC, so-called second generation, can be fed by it as in the Vattenfall case in Berlin where 120MW have been recently installed. The detailed characterization of the units is done with a piecewise linearization of their useful effect as a function of one and two variables, with their start-up costs and ramp constraints. Electric transmission constraints are also included. The large-scale thermal energy storage follows practices employed in modern DH systems [11].

## 2 Mathematical model

In this section, a mathematical model for heat and power generation unit commitment co-optimizing the electric transmission system and its integration with thermal areas is presented.

### 2.1 Objective function

The objective of the optimization model is to determine the minimum operational cost of the joint energy systems for the day-ahead unit commitment, $t \in \{1, \dots 24\}$, given by:

$$\min. \sum_t (c_t^f + c_t^{st} + c_t^{OM}). \quad (1)$$

The system's operational cost can be divided into fuel consumption $c_t^f$ (2), start-up costs $c_t^{st}$ (3), and operation and maintenance costs $c_t^{OM}$ (4):

$$c_t^f = \sum_i \hat{C}_i^f f_{i,t}, \quad \forall t \quad (2)$$

$$c_t^{st} = \sum_i \hat{C}_i^{st} \tau_{i,t}, \quad \forall t \quad (3)$$

$$c_t^{OM} = \sum_i (\hat{C}_i^{OM,tm} \theta_{i,t} + \hat{C}_i^{OM,st} \tau_{i,t} + \hat{C}_i^{OM,f} f_{i,t}), \quad \forall t. \quad (4)$$

### 2.2 Energy systems modelling

#### 2.2.1 Electric energy system

The nodal power balance based on a DC power flow is presented in (5). Expression (6) represents the power flow through a transmission line. The static thermal rating (STR) is the maximum permissible power through a line (7).

$$P_{n,t} - \sum_{i \in EP} \Omega_{i,n} p_{i,t} + \sum_{w} p_{n,w,t} + W_n - W_n^{spill} = 0,$$
$$\forall n,t \quad (5)$$
$$B_{n,w}(\delta_{n,t} - \delta_{w,t}) = p_{n,w,t}, \quad \forall n,w,t \quad (6)$$
$$-\overline{P}_{n,w} \leq p_{n,w,t} \leq \overline{P}_{n,w}, \quad \forall n,w,t. \quad (7)$$

### 2.2.2 Thermal energy system

The thermal energy balance for each zone $z$ of the system is given by (8). The total thermal load per zone is represented by $H_{z,t}^d$. The energy level of the zonal thermal storage is given by $s_{HT,z,t}$. $\Gamma_{i,z}$ is a binary parameter representing the zone served by unit $i \in HT$. The parameter $\eta_z$ represents the heat transfer efficiency within zone $z$, and $L_z^s$ the storage's hourly energy loss.

$$H_{z,t}^d = \eta_z \left( \sum_{i \in HT} \Gamma_{i,z} h_{i,t} + [s_{z,t} - s_{z,t+1}] \right) - L_z^s s_{z,t},$$
$$\forall z,t. \quad (8)$$

## 2.3 Units characterization

### 2.3.1 Start-up procedures

Equations (9)-(11) define the value of the binary variable $\tau_{i,t}$, representing whether unit $i$ undergoes a start-up procedure during $t$, thereby ensuring that it equals to 1 only for the time step when the unit $i$ is turned on. The binary variable $\theta_{i,t}$ states a unit's on/off status.

$$\tau_{i,t} \leq \theta_{i,t}, \quad \forall i,t \quad (9)$$
$$\tau_{i,t} \leq 1 - \theta_{i,t-1}, \quad \forall i,t \quad (10)$$
$$\tau_{i,t} \geq \theta_{i,t} - \theta_{i,t-1}, \quad \forall i,t. \quad (11)$$

The number of start-up and shut-down procedures for such large thermal units needs to be limited for reliability. The maximum number of start-up procedures during a day is set by:

$$\sum_t \tau_{i,t} \leq N_i^{st}, \quad \forall i,t \quad (12)$$

The minimum number of consecutive commitment periods for a unit is given by:

$$\theta_{i,t} \geq \sum_{\tilde{t}=1}^{MT_i} \tau_{i,t-\tilde{t}}, \quad \forall i,t \quad (13)$$

### 2.3.2 Performance curves for units with one independent variable

In most generation units, a nonlinear dependency exists between fuel consumption and energy production. To model the performance curves as convex functions, they are characterized by a piecewise linear approximation (14)-(19) [12]. The performance curves are sampled by $J_i$ break-points, and are modelled via the convex combination of the sampling points. Fuel consumption is given by (17). Constraints (14)-(18) apply to units producing electric energy, while (14)-(17) and (19) model those producing heat. CHP units are modeled by (14)-(19). P2H units can be characterized by setting the negative values to the electricity generation's sampling points.

$$\alpha_{i,j,t} \leq \beta_{i,j,t-1} + \beta_{i,j,t}, \quad \forall i \in I_1, j, t \quad (14)$$
$$1 = \sum_{j=1}^{|J_i|-1} \beta_{i,j,t}, \quad \forall i \in I_1, j, t \quad (15)$$
$$1 = \sum_j \alpha_{i,j,t}, \quad \forall i \in I_1, j, t \quad (16)$$
$$f_{i,t} = \sum_j \alpha_{i,j,t} \hat{F}_{i,j}, \quad \forall i \in I_1, t \quad (17)$$
$$p_{i,t} = \sum_j \alpha_{i,j,t} \hat{P}_{i,j}, \quad \forall i \in EP \cap I_1, t \quad (18)$$
$$h_{i,t} = \sum_j \alpha_{i,j,t} \hat{H}_{i,j}, \quad \forall i \in HT \cap I_1, t. \quad (19)$$

### 2.3.3 Performance curves of units with two independent variables

Units with two degrees of freedom are those whose useful effect depends on the value of two independent variables, $f_{i,t}$ and $o_{i,t}$. The characteristic curves are sampled via $J_i$ and $K_i$ break-points for the independent variables [12]. Constraints (20)-(26) model the units' performance curves:

$$\alpha_{i,j,k,t} \leq \beta_{i,j-1,k-1,t}^{up} + \beta_{i,j-1,k-1,t}^{low} + \beta_{i,j,k-1,t}^{up} + \beta_{i,j-1,k,t}^{low}$$
$$+ \beta_{i,j,k,t}^{up} + \beta_{i,j,k,t}^{low}, \quad \forall i \in I_2, j, k, t \quad (20)$$
$$1 = \sum_{j=1}^{|J_i|-1} \sum_{k=1}^{|K_i|-1} \left( \beta_{i,j,k,t}^{up} + \beta_{i,j,k,t}^{low} \right), \quad \forall i \in I_2, t \quad (21)$$
$$1 = \sum_j \sum_k \alpha_{i,j,k,t}, \quad \forall i \in I_2, t \quad (22)$$
$$f_{i,t} = \sum_j \sum_k \alpha_{i,j,k,t} \hat{F}_{i,j}, \quad \forall i \in I_2, t \quad (23)$$
$$o_{i,t} = \sum_j \sum_k \alpha_{i,j,k,t} \hat{O}_{i,k}, \quad \forall i \in I_2, t \quad (24)$$
$$p_{i,t} = \sum_j \sum_k \alpha_{i,j,k,t} \hat{P}_{i,j,k}, \quad \forall i \in EP \cap I_2, t \quad (25)$$
$$h_{i,t} = \sum_j \sum_k \alpha_{i,j,k,t} \hat{H}_{i,j,k}, \quad \forall i \in HT \cap I_2, t. \quad (26)$$

### 2.3.4 Technical limits for generation and storage

The electric and thermal generation limits are respectively presented in (27) and (28). For the thermal storage, the heat accumulated must be non-negative and below its technical capacity, as shown in (29). Ramp rates limit the power change in consecutive periods, represented by a ramp-up rate - $UR_i$, and ramp-down rate - $DR_i$; modelled by (30).

$$\theta_{i,t} \underline{E_i} \leq e_{i,t} \leq \theta_{i,t} \overline{E_i}, \quad \forall i \in EP, t \quad (27)$$

$$\theta_{i,t}\underline{H_i} \leq h_{i,t} \leq \theta_{i,t}\overline{H_i}, \quad \forall i \in HT, t \quad (28)$$

$$0 \leq s_{z,t} \leq \overline{S_z}, \quad \forall z, t \quad (29)$$

$$-DR_i \leq e_{i,t} - e_{i,t-1} \leq UR_i, \quad \forall i \in EP, t \quad (30)$$

## 3 Numerical Tests

To demonstrate the impact of a joint electric and thermal system operation on transmission decongestion, the day-ahead unit commitment of an electric transmission network and two heating zones is presented. The electric system is a modified IEEE 24-bus test system [13]. The two heating zones represent large DH networks fed by one CHP plant each. CHP-2d and CHP-1d have been added to the electric system at nodes 15 and 23, and thermal zones 1 and 2, respectively. The CHP-1d is a gas turbine with heat recovery, while the CHP-2d is a natural gas combined cycle with an extraction condensing steam turbine, i.e., its two independent variables are fuel consumption - $f_{i,t}$, and valve opening - $o_{i,t}$. The two thermal zones representing the heating system are located at two nodes of the electrical system with a given load profile. Each zone accounts two gas fired boilers, one CHP unit, one electric boiler (P2H) and a zonal thermal energy storage. The thermal capacity of the CHPs is taken to be 45 % of the maximum load and the gas fired boilers act as a back-up to fulfil the entire load according to the operational practice of district heating networks [11]. The system's technical parameters are implemented as per [14] and can be found in [15]. The simulations were performed using the scientific computing language Julia v1.0.5., with the optimization package JuMP v0.21.2 and Cbc v0.6.2 as a MILP solver, on a computer with an Intel Xenon Gold 6148 CPU @ 2.4GHz and 256 GB RAM.

### 3.1 Renewable integration index

In addition to the operational cost, a performance index is used to assess the system operation performance under congestion, namely the system utilization (SU) [16]. SU is the ratio between the total system generation in the congested and uncongested cases. Given that our model is based on the use of renewable generation in a lossless transmission system, the use of the renewable integration (RI) index is proposed that compares the share of renewable energy employed in the constrained and unconstrained cases:

$$RI = \frac{\sum_n \left(W_n - W_n^{\text{spillage}}\right)}{\sum_n \left(W_n - W_n^{\text{spillage,unconstrained}}\right)} - 1 \quad (31)$$

The RI indicates how much more wind energy is being integrated with the use of flexibility measures, when compared to the unconstrained decoupled case with the same wind penetration level. A negative RI indicates less wind integration at the optimal operation.

### 3.2 Results and discussion

The numerical tests have been divided into seven cases to highlight the effects on renewable energy usage and transmission congestion, increasing integration of the electric with thermal infrastructure:

- *Case 0:* the electric and thermal systems are dispatched in a decoupled manner. First, the operation of the heating zones is optimized. With the thermal dispatch, the CHP units' thermal generation is fixed, and the unit commitment of the electric system is performed without transmission constraints. The maximum wind penetration is set at 45% for cases 0 to 4.
- *Case 1:* this case presents the current operational practice, where the electric and thermal systems are decoupled as in *Case 0*, with transmission constraints in the electric network.
- *Case 2:* the operation scheduling of the thermal and electric systems is co-optimized, as in Section 2.
- *Case 3:* the system is optimized as in *Case 2*. Thermal energy storage is added to the heat zones.
- *Case 4:* electric boilers are included in the joint operation presented in *Case 3*.
- *Case 5:* same as 4 with wind penetration of 65%.
  *Case 6:* same as 4 with wind penetration of 85%.

The optimization results are presented in **Table 1**. In the decoupled case without transmission constraints, *Case 0*, three lines of the IEEE 24-bus reliability test system exceed their rated capacity with an average flow of 119% of their STR, see **Fig. 1**. The maximum line usage is presented at 133% for line *15-21*.

**Table 1.** Main model features and optimization results

|  | Case 0 | Case 1 | Case 2 | Case 3 | Case 4 | Case 5 | Case 6 |
|---|---|---|---|---|---|---|---|
| **Wind penetration** | 45% | 45% | 45% | 45% | 45% | 65% | 85% |
| **Transmission constraints** |  | X | X | X | X | X | X |
| **Coupled systems** |  |  | X | X | X | X | X |
| **Storage** |  |  |  | X | X | X | X |
| **Power-to-Heat** |  |  |  |  | X | X | X |
| **Binary variables** | 809 |  | 779 | 841 | 889 |  |  |
| **Total number of variables** | 18 024 |  | 17 424 | 17 520 | 17 856 |  |  |
| **Computational time [s]** | 6.98 | 11.79 | 12.11 | 1 264.31 | 1 975.47 | 372.20 | 166.11 |
| **Total cost [million €]** | 4.407 | 4.758 | 4.644 | 4.295 | 4.294 | 2.417 | 1.245 |
| **Cost variation [%]** | -7.37 | — | -2.40 | -10.49 | -10.5 | -49.20 | -73.83 |
| **Fuel consumption [MWh]** | 60 516.4 | 69 019.4 | 66 248.6 | 61 287.42 | 61 269.0 | 34 469.3 | 17 748.9 |
| **Wind spillage [MWh]** | 227.3 | 227.3 | 0.0 | 20.5 | 0.00 | 450.9 | 1 218.5 |
| **Renewable integration (RI) [%]** | — | 0 | 1.0 | 0.9 | 1.0 | 7.5 | 30.6 |

Even though *Case 0* has no limits in transmission capacity, there exists wind spillage, as seen in **Fig. 2**. The wind spillage is a consequence of the minimum technical limits of the generation and the start-up costs. In this case, generator eight must generate at minimum capacity during hours 6 and 7, avoiding additional start-up costs, thereby leading to wind spillage in those hours.

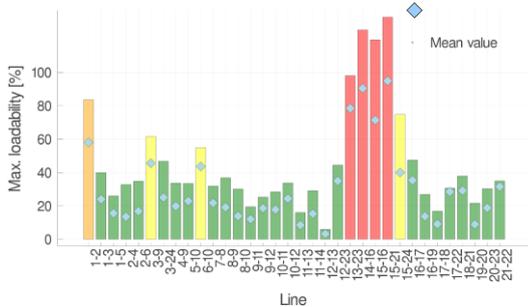

**Fig. 1.** Transmission lines usage for Case 0

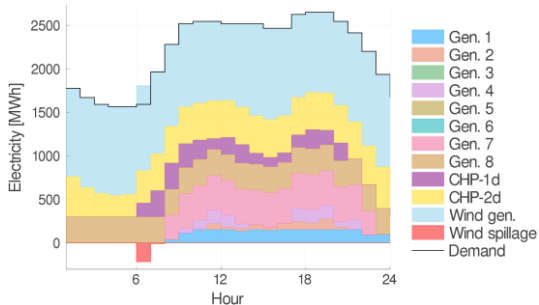

**Fig. 2.** Electric generation scheduling for Case 0

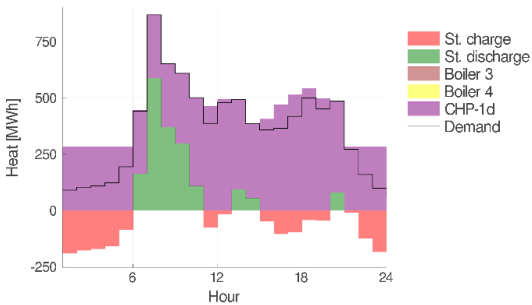

**Fig. 3.** Thermal dispatch of zone 2 in *Case 3*

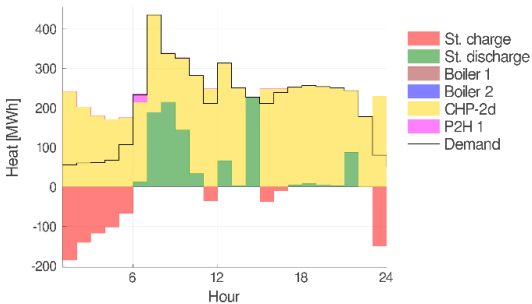

**Fig. 4.** Thermal dispatch in zone 1

The decoupled and constrained case, *Case 1*, has a 0% RI, as seen in **Table 1**. Even though this case introduces transmission constraints, the system has enough transmission capacity to integrate the same wind energy amount without exceeding flow limits. However, since the electric generation is dispatched differently to ensure adequate line flow, its operating cost increases by 7.37% compared to *Case 0*.

When the operation of electric and thermal systems is coupled, *Case 2*, the flexible use of the CHP-2d, allows complete wind energy use. The total wind spillage amounts to 0 MWh, yielding an RI of 1% and reducing the total operation cost by 2.4%.

In *Case 2,* the CHP-1d is not used since it cannot be scheduled so that the cost reductions derived from additional electricity generation and consequent extra heat can compete with the gas boilers. This issue is overcome once the thermal energy storage is introduced in *Case 3*. As seen in **Fig. 3**, the active use of the storage system in zone 2 allows storage of the excess heat generated by the CHP-1d during off-peak hours and then releases it in peak hours to complement the unit's capacity. Thus, the thermal energy storage enhances the CHP-1d efficiency above that of the combination of gas boilers and electric generators. The active use of the thermal energy storage further reduces the total costs by 10.49% compared to *Case 1*. However, since the CHP-1d is employed, 20.5 MWh of the wind generation must be spilled, leading to a RI of 0.9%.

The introduction of P2H, in the form of electric boilers, results in zero wind spillage. However, as seen in **Fig. 4**, they are only used in zone 1, and its consumption equals the wind spillage of *Case 3*. The P2H is used in zone 1 since the CHP-2d, unlike the CHP-1d in zone 2, can reduce its heat generation for the sixth hour without affecting its electricity production.

Due to its low coefficient of performance (COP) of 0.99, the electric boiler is only used in *Case 4* to prevent the wind spillage and not to replace the heat generated from the CHPs actively. Such P2H's usage reflects a wind penetration not high enough to make the price of electricity so low to gain vs. CHP units at low load. Once the wind penetration is increased to 65%, Case 5 shows an active complement between P2H and CHP-2d unit in zone 1, **Fig. 5**. However, the P2H is only used during night hours when there is low electric demand, **Fig. 6**, and the CHP-2d is not needed to support the electricity generation as it is not convenient to sell electricity.

In *Case 5*, the P2H unit in zone 2 is not utilized since the CHP-1d, even with the storage, is not flexible enough to actively manage its power-to-heat ratio. When the wind penetration is increased to 85% in *Case 6*, the high amount of wind energy available makes it economically feasible to stop using the CHP-1d and replace its generation with the electric and gas boilers, in combination with the energy storage unit, see **Fig. 7**.

## Conclusions

In this work, a mixed-integer linear programming (MILP) model is proposed for optimal scheduling of electric transmission systems integrated with heating systems, namely zones where large district heating (DH) is present. The model characterizes the performance of electricity- and heat-only units, cogeneration units with one and two degrees-of-freedom, with the latter allowing to decouple in a flexible fashion thermal and electric generation, large scale power-to-heat (P2H) units, and thermal energy storage.

Seven case studies are assessed, and the optimal system scheduling under different operation conditions is presented for each of them. The co-optimization of the electric and thermal systems reduced the operational costs by 2.4% compared to the decoupled operation scheduling; furthermore, it allows integrating renewable generation without wind spillage. Higher cost savings of 10.49% are obtained using thermal energy storage, which improves the efficiency of the otherwise inflexible one degree-of-freedom (d.o.f.) combined heat and power (CHP) unit, enabling to decouple its thermal and electric generation. Finally, large scale power-to-heat units above 20 MW-electric, deployable at the transmission scale are considered. To be realistic, electric boilers are selected, showing that they are not profitable for a low wind penetration of 45% and becoming a viable option for higher penetrations of 65 and 85%. P2H, for retrofitting existing DH networks (above 100ºC), still needs technological development, e.g., industrial heat pumps, else P2H electric boilers are viable only for high renewable penetration and during low electric demand hours. For low wind penetration, one-d.o.f. units complemented with heat storage and two-d.o.f. systems showed the best results.

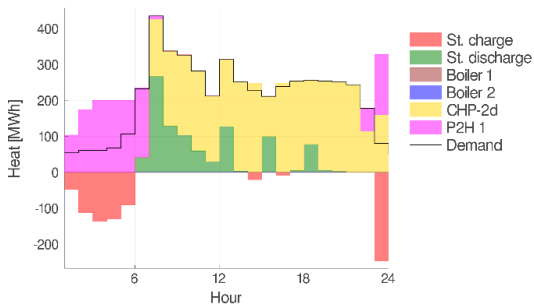

**Fig. 5.** Thermal dispatch in zone 1 on *Case 5*

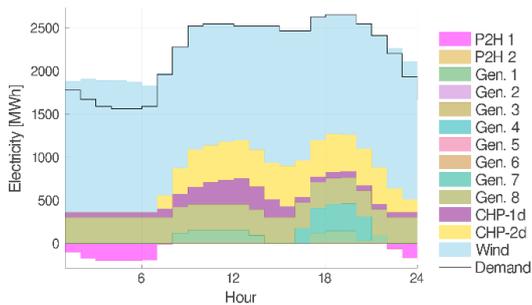

**Fig. 6.** Electric dispatch on *Case 5*

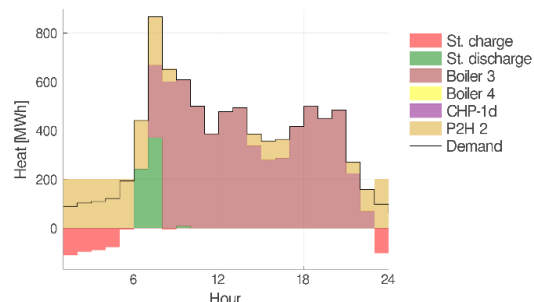

**Fig. 7.** Thermal dispatch in zone 2 on *Case 6*